\documentclass{amsart}
 \newtheorem{thm}{Theorem}[section]
 
 \newtheorem{lem}[thm]{Lemma}
 
 \theoremstyle{definition}
 \newtheorem{defn}[thm]{Definition}
 \theoremstyle{remark}
 \newtheorem{rem}[thm]{Remark}
 
 \numberwithin{equation}{section}

\begin{document}

\title[Riesz and Rajchman sets]{Some supports of Fourier transforms of \\ singular measures are not Rajchman}%
\author{Maria Roginskaya}%
\address{Mathematical Sciences\\
Chalmers University of Technology\\
SE-41296, Gothenburg, Sweden\\
 \\
Mathematical Sciences\\
Gothenburg University\\
SE-41296, Gothenburg, Sweden}%
\email{maria.roginskaya@chalmers.se}%

\thanks{This work was accomplished with the support of Fondation Sciences Math\'ematiques de Paris}%
\keywords{Rajchman sets, Riesz sets, Riesz products, singular measures, support of Fourier transform}%

%\date{}%
%\dedicatory{}%
%\commby{}%
% ----------------------------------------------------------------
\begin{abstract}
The notion of Riesz sets tells us that a support of Fourier transform of a measure with non-trivial singular part has to be large. The notion of Rajchman sets tells us that if the Fourier transform tends to zero at infinity outside a small set, then it tends to zero even on the small set. Here we present a new angle of an old question: Whether every Rajchman set should be Riesz.
\end{abstract}
\maketitle
% ----------------------------------------------------------------
\section{Introduction}

The consideration of the properties of measures and their Fourier transforms is a classical area of Harmonic Analysis. In particular the following is well known.

\begin{thm}[Rajchman, 1929 \cite{Raj}] If for a finite measure $\mu$ on the unit circle $\Bbb T$ holds $\widehat{\mu}(n)\rightarrow 0$ when $n\rightarrow -\infty$, then it holds also that $\widehat{\mu}(n)\rightarrow 0$ when $n\rightarrow +\infty$.
\end{thm}

This motivates the following.
\begin{defn} We say that $\Lambda\subset {\Bbb Z}$ is a {\it Rajchman set} if as soon as $\widehat{\mu}(n)\rightarrow 0$
when $|n|\rightarrow +\infty, n\in {\Bbb Z}\setminus \Lambda$, then $\widehat{\mu}(n)\rightarrow 0$ when $|n|\rightarrow +\infty, n\in \Lambda$.
\end{defn}

With this definition the Rajchman theorem says that the non-negative integers is a Rajchman set.

Now, given a (signed) Radon measure $\mu$ on the unit circle $\Bbb T$, we can present it as $\mu=f \cdot m + \mu_s$, where $m$ is the Lebesgue measure and $\mu_s$ is the singular with respect to Lebesgue measure part of the measure $\mu$. We known the following.

\begin{thm}[F. and M. Riesz's, 1916, \cite{Ri}] If a finite measure $\mu$ has the property $\widehat{\mu}(-n)=0$ for $n=1,\ldots$, then the measure is absolutely continuous with respect to Lebesgue measure, i.e. $\mu=f \cdot m$, where $f\in L^1({\Bbb T})$.
\end{thm}

This result motivates the following definition.

\begin{defn} We say that a subset $\Lambda \subset {\Bbb Z}$ is a Riesz set if it has the property, that if $supp(\widehat{\mu})\subset\Lambda$ then $\mu$ has no singular part.
\end{defn}

With this definition the F. and M. Riesz theorem says that the non-negative integers is a Riesz set.

\begin{thm}[Host, Parreau, 1978 \cite{HP}\footnote{It is actually proven in \cite{HP} not only for $\Bbb T$ but for any compact group}.
] A set $\Lambda\subset {\Bbb Z}$ is a Rajchman set iff it doesn't contain any shift of the Fourier support of a Riesz product, i.e. any set $\Omega((n_j))=\{\sum \epsilon_j n_j: \epsilon_j=-1,0,1; \sum |\epsilon_j|<\infty\}$, where $(n_j)$ is an infinite sequence.
\end{thm}
Thus, any set which is not Rajchman, contains the support of the Fourier transform of a singular measure, and thus is not Riesz (or, without negations, that every Riesz set is a Rajchman set).

A natural question is following: Is every Rajchman set a Riesz set? (I.e. Do the classes of Riesz and Rajchman sets coincide?) As far to the author's knowledge, this question was first raised by Pigno, 1978 \cite{P}.

As we are unable to answer the question, we want to diversify it:

\begin{defn} We say that a closed set $E\subset{\mathbb T}$ is a {\it parisian set} if for every non absolutely continuous measure $\mu\in M(E)$, the support of it's Fourier transform is not a Rajchman set.
\end{defn}

The original question thus becomes: Is $\mathbb T$ a parisian set?

While we are not able to answer the question above, we can show that some parisian sets do exist. As any subset of a parisian set is parisian, it is clear that a positive answer on the original question would imply all the results we prove here. Yet, there are good chances that the answer is negative and a negative answer would give the study of the parisian sets some interest.

It is natural to expect that the parisian sets should be "small". Thus we try to construct a "big" parisian set.

{\bf Main Theorem A.} {\it For any $\alpha<1$ there exists a closed parisian set $E$, such that $dim_H(E)\geq\alpha$, where $dim_H(E)$ means the Hausdorff dimension of $E$.}

{\bf Main Theorem B.} {\it For any $\alpha<1$ there exists a Borel parisian set $E$ such that it is an additive subgroup of $\Bbb T$ and $dim_H(E)\geq\alpha$.}

\subsubsection*{Notations.} In what follows we identify $\Bbb T$ with $(-1,1]$, so that the Fourier coefficients are $\widehat{\mu}(n)=\frac12\int e^{i\pi nx}d\mu(x)$.
\section{Construction of a big parisian set}

Let us first introduce a test to establish that a set is parisian.
\begin{lem}\label{test} If there exist $\delta>0$ and a sequence $(N_j)_{j=1}^\infty$ such that for every $j$ the set $E$ is a subset of $\frac2{N_j}{\Bbb Z}+[-1/2N_j^{1+\delta},1/2N_j^{1+\delta}]$, then the set $E$ is parisian.
\end{lem}
\begin{proof} Let us fix $\mu\in M_s(E)$. We want to show that $supp(\widehat{\mu})$ contains a shift of a set $\Omega((n_j))$. Up to a shift of the Fourier transform we may assume without loss of generality that $\widehat{\mu}(0)\neq 0$.

Here we construct the sequence $(n_j)$ as a subsequence of $(N_j)$ inductively. Assume that $(k-1)$ first terms of the sequence $(n_j)$ are chosen. This means that for all combinations of $\epsilon_j$ the sum $\sum\limits_0^{k-1}\epsilon_j n_j\in supp(\widehat \mu)$. Thus, we know that $\int  e^{i\pi \sum\limits_{j=1}^{k-1}\epsilon_jn_j x}d\mu(x)\neq 0$, for all combinations $(\epsilon_j=-1,0,1)_{j=1}^{k-1}$. We can take $\gamma_{k-1}$ to be the minimum of the absolute value of the $3^{k-1}$ non-zero numbers, so that $|\int  e^{i\pi \sum\limits_{j=1}^{k-1}\epsilon_jn_j x}d\mu(x)|\geq \gamma_{k-1}$. We want to show that for some sufficiently large $n_k=N_{j_k}$ for all combinations of $\epsilon_j$ holds $\int  e^{i\pi \sum\limits_{j=1}^{k}\epsilon_jn_j x}d\mu(x)\neq 0$.

Indeed, as $E\subset 2{\mathbb Z}/N_m+[-1/N_m^{1+\delta},1/N_m^{1+\delta}]$, we know that $|e^{i\pi (\pm N_m x)}-1|\leq \frac{\pi}{N_m^\delta}$, when $x\in E$. Now we see that $$|\int_E  e^{i\pi \sum\limits_{j=1}^{k}\varepsilon_jn_j x}d\mu(x)-\int_E  e^{i\pi \sum\limits_{j=1}^{k-1}\varepsilon_jn_j x}d\mu(x)|\leq \int_E |d\mu||e^{i\pi \pm N_m x}-1|\leq \|\mu\|\frac{\i}{N_m^\delta}.$$ Thus, for sufficiently large $m$ we can be sure that the later is less than $\frac12\gamma_{k-1}$. Now, we see that by the triangle inequality $|\int_E  e^{i\pi \sum\limits_{j=1}^{k}\varepsilon_jn_j x}d\mu(x)|\geq \frac12 \gamma_{k-1}>0$ for all the combinations of $\epsilon_j=-1,0,1$, with $j=1,\ldots, k$, and $n_k=N_m$.
\end{proof}

A slight modification of the proof gives us the following.
\begin{lem} For an increasing sequence $(N_j)\subset {\Bbb N}$ and $\delta>0$ the set $\widetilde{E}=\{x\in {\Bbb T}: \sup\limits_j(dist(x,2{\Bbb Z}/N_j)/N_j^{1+\delta})<\infty\}$ is a parisian set.
\end{lem}
\begin{proof} We start from observing that $\widetilde{E}=\bigcup\limits_{t\in {\Bbb N}}E_t$, where $$E_t=\{x\in {\Bbb T}: \sup\limits_j(dist(x,2{\Bbb Z}/N_j)/N_j^{1+\delta})\leq t\}$$ is an increasing sequence of closed sets.

Now, we start the proof exactly as the previous one, but after the choice of $\gamma_{k-1}$ and before the choice of $n_k$ we do one more step: We pick $t_k$ large enough that $\mu_k=\mu|_{E_{k}}$ satisfies $\|\mu-\mu_k\|<\frac13\gamma_{k-1}$. Then we see that $|\int  e^{i\pi \sum\limits_{j=1}^{k-1}\epsilon_jn_j x}d\mu_k(x)|\geq \frac23\gamma_{k-1}$. We proceed in the same way as before with $\mu_k$ in place of $\mu$, and find $n_k=N_{m_k}$ such that $|\int_E  e^{i\pi \sum\limits_{j=1}^{k}\varepsilon_jn_j x}d\mu_k(x)|\geq \frac12 \gamma_{k-1}$. Then, $|\int_E  e^{i\pi \sum\limits_{j=1}^{k}\varepsilon_jn_j x}d\mu(x)|\geq \frac16 \gamma_{k-1}>0$.
\end{proof}

\begin{rem} The set $\widetilde{E}$ is obviously an additive subgroup of $\Bbb T$ and thus either finite or dense in $\Bbb T$.
\end{rem}

Let us now construct a set $E$ of large Hausdorff dimension which satisfies the hypothesis of the Lemma \ref{test}, and is thus parisian. As the constructed set is a subset of $\widetilde{E}$ it will also give us the estimate\footnote{This estimate is well known, but we give the proof for the sake of completeness.} on the Hausdorff dimension of $\widetilde{E}$. Fix $\alpha\in (0,1)$, and choose $\delta>0$ so that $\delta=1-\alpha$. We will construct a rapidly increasing sequence $\{N_j\}$, and related sequence of closed sets $C_j\subset (-1,1)$,  such that the sets $C_j$ is the union of the closed intervals with centrums in $2{\Bbb Z}/N_j$, of length $1/N_j^{1+\delta}$ which are entirely contained in $\bigcap\limits_{k=1}^{j-1}C_{k}$. We will let then the set $E=\bigcap_j C_j$, which is obviously closed. The set constructed in such a way is a Cantor-type set, and we show that provided the sequence $N_j$ grows quickly enough the dimension of such a set is at least $\alpha$.

\begin{lem} $dim_H(E)\geq \alpha$.
\end{lem}
\begin{proof} In order to prove that the Hausdorff dimension of $E$ is at least $\alpha$ we will show that it is at least $s$ for any $0<s<\alpha$, and to do so we construct a finite measure $\mu$ supported on $E$ such that $\mu(I)\leq c_s |I|^s$ for any interval $I$ (it is  a standard fact of Geometric Measure Theory that a measure satisfying such an estimate should have support of Hausdorff dimension at least $s$, see for example \cite{Matt}).

Let us take a subset $D_k$ of $\bigcap\limits_{j=1}^k C_j$, which is a collection of intervals of length $1/N_k^{1+\delta}$. This collection is defined inductively: we know that every interval of length $1/N_{k-1}^{1+\delta}$ contains at least $N_k/2N_{k-1}^{1+\delta} -1$ points of $2{\mathbb Z}/N_k$. Thus, every interval of $D_{k-1}$ contains (entirely) at least $M_k=N_k/2N_{k-1}^{1+\delta} - 3$ intervals with centrum in $2{\mathbb Z}/N_k$ and length $1/N_k^{1+\delta}$. (To make the estimates more simple we assume $(N_k)$ to grow so rapidly that $M_k\geq N_k/4N_{k-1}^{1+\delta}$.)

We pick from each interval of $D_{k-1}$ exactly $M_k$ such intervals. All together we will have picked $M_k\prod\limits_{j=1}^{k-1}M_j$ intervals of length $\frac1{N_{k}^{1+\delta}}$. Then we take the probability measure $\mu_k$ equally distributed on the $\prod\limits_{j=1}^k M_j$ intervals of $D_k$. We introduce $\mu$ as a weak limit point of $\mu_k$ (which has to be a probability measure supported by $E=\cap C_j$).

Let us estimate $\mu(I)$ where $1/N_{k-1} >|I|\geq 1/N_k$. The interval can intersect at most $N_k |I|/2 +3$ intervals of $D_k$ (as $N_k|I|\geq 1$, we may use that it is at most $4 N_k|I|$ intervals). As the measure of each interval of $D_k$ is $1/\prod^k_{j=1} M_j$ we see that $\mu(I)\leq 4|I|N_k /\prod_{j=1}^k M_k$ where $$\prod M_k\geq (N_k/N_1)/(4^{k-1} (\prod_{j=1}^{k-1}N_j)^\delta).$$ Thus, $\mu(I)\leq N_1 4^k (\prod_{j=1}^{k-1} N_j)^\delta |I|=N_1 4^k (\prod_{j=1}^{k-1}N_k)^\delta |I|^{1-s} |I|^s$.

Our task is fulfilled if we show that $c_{k,s}=N_1 4^k(\prod_{j=1}^{k-1}N_j)^{\delta}|I|^{1-s}$ is bounded above independently from $k$. We know that $|I|<1/N_{k-1}$, and, as $\delta=1-\alpha$, we see that $c_{k,s}\leq N_1 4^k (\prod_{j=1}^{k-2} N_j)^\delta / N_{k-1}^{\alpha-s}$. It remains to take the sequence $(N_k)$ such that $(N_1 4^{k+2}(\prod_{j=1}^k N_j)^\delta)^k< N_{k+1}$. For any fixed $s$ the sequence $c_{k,s}$ tends to zero, and so is bounded. (Notice that the bound $c_s=\sup\limits_k \{c_{k,s}\}$ grows as $s\rightarrow \alpha$, but we only need it to be finite.)
\end{proof}

% ----------------------------------------------------------------

\end{document}